\documentclass[10pt]{article}
\usepackage[utf8]{inputenc}
\usepackage[slovak,english]{babel}
\usepackage{a4}
\usepackage{subfig}
\usepackage{graphicx}
\usepackage{cite}
\usepackage{caption}
\usepackage{color}
\usepackage{amsmath, amsthm, amssymb}
\usepackage{mathtools}
\usepackage{multicol}
\usepackage{url}
\usepackage{setspace}
\usepackage{pdfpages}
\usepackage{afterpage}
\usepackage{hyperref}
\usepackage[total={17cm,25cm}, top=25mm, bottom=25mm, right=25mm, left=40mm, includefoot]{geometry}

\begin{document}
\title{Interactive 4-D Visualization of Stereographic Images From the Double Orthogonal Projection.}

\author{Michal~Zamboj\\ 
Charles University, Faculty of Education\\
M. D. Rettigové 4, 116 39 Prague 1.}
\date{}
\maketitle
\begin{abstract}
The double orthogonal projection of the 4-space onto two mutually perpendicular 3-spaces is a method of visualization of four-dimensional objects in a three-dimensional space. We present an interactive animation of the stereographic projection of a hyperspherical hexahedron on a 3-sphere embedded in the 4-space. Described are synthetic constructions of stereographic images of a point, hyperspherical tetrahedron, and 2-sphere on a 3-sphere from their double orthogonal projections. Consequently, the double-orthogonal projection of a freehand curve on a 3-sphere is created inversely from its stereographic image. Furthermore, we show an application to a synthetic construction of a spherical inversion and visualizations of double orthogonal projections and stereographic images of Hopf tori on a 3-sphere generated from Clelia curves on a 2-sphere.
\end{abstract}

\section{Introduction}
The double orthogonal projection is a four-dimensional generalization of Monge's projection. In Monge's projection, an object in the 3-space is orthogonally projected into two mutually perpendicular planes (horizontal and vertical, i.e. top and front view). One of these planes is chosen to be the drawing (or picture) plane and the second is rotated about their intersecting line to the drawing plane. Therefore, each point in the 3-space has two conjugated images in the drawing plane. In the double orthogonal projection, an object is in the 4-space, and we project it orthogonally into two mutually perpendicular 3-spaces. Let $x,y,z$, and $w$ be the orthogonal system of coordinate axes of the 4-space, and $\Xi(x,y,z)$ and $\Omega(x,y,w)$ be the 3-spaces of projection. If $\Omega$ is chosen to be the modeling 3-space (instead of the drawing plane), then $\Xi$ is rotated about their common plane $\pi(x,y)$ such that $z$ and $w$ have opposite orientations (usually $w$ up and $z$ down). Analogically, each point in the 4-space has two conjugated images in the modeling 3-space. Elementary constructions and principles of the double orthogonal projection were described in \cite{Zamboj2018a}, sections, and lighting of polytopes in \cite{Zamboj2018b}, intersections of lines, planes, and 3-spaces with a 3-sphere in \cite{Zamboj2019b}, and regular quadric sections as intersections of 4-dimensional cones with 3-spaces in \cite{Zamboj2019a}. Similarly to Monge's projection, in which an observer can reach any point of the drawing plane, in the four-dimensional case an observer can reach any point of the modeling 3-space. For this purpose, our constructions are supplemented by online interactive models \cite{Zamboj2020icggGGB} created in \emph{GeoGebra~5}. 

\begin{figure}[!htb]
\centering
\includegraphics[height=4.5cm]{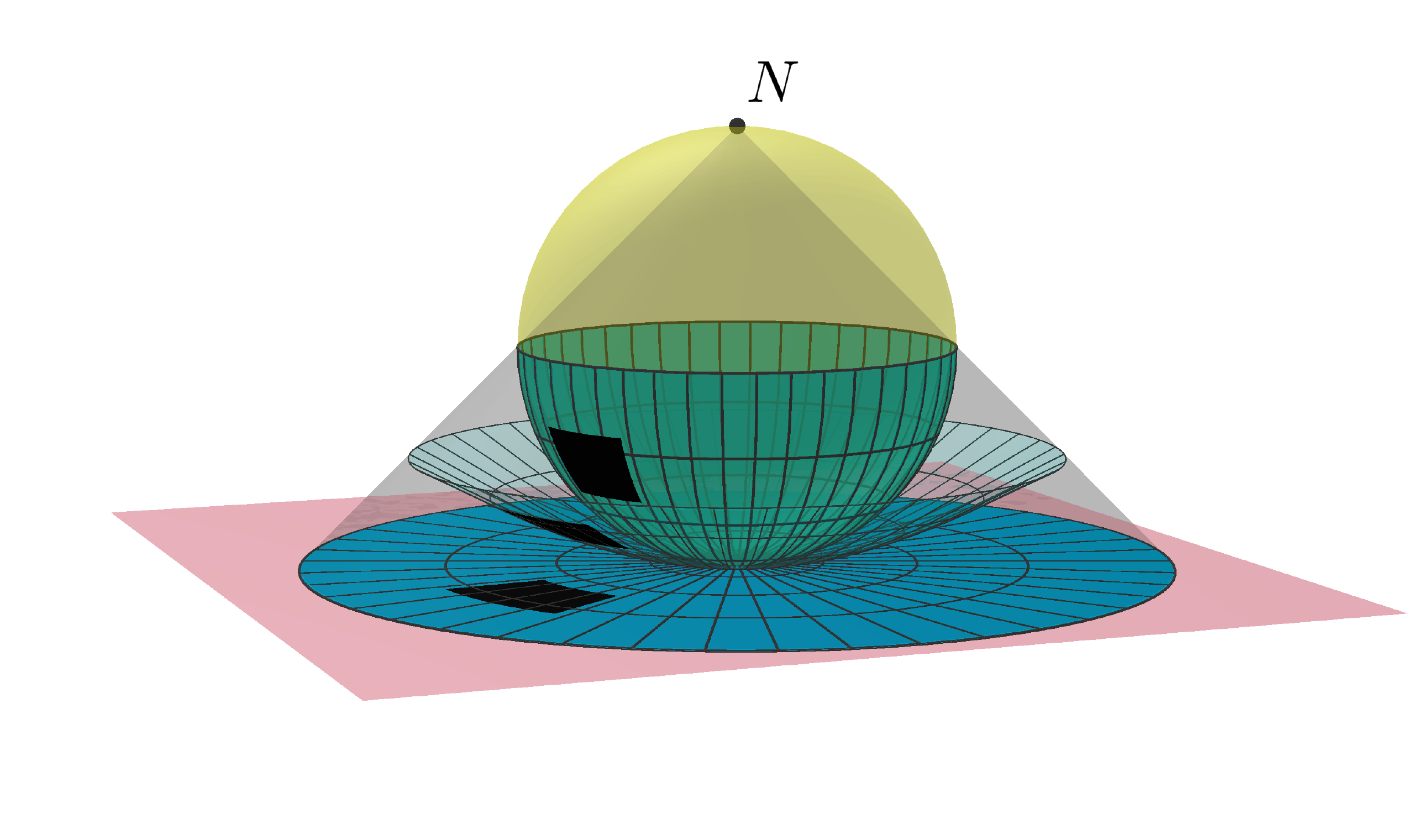}\hfill\includegraphics[height=4.5cm]{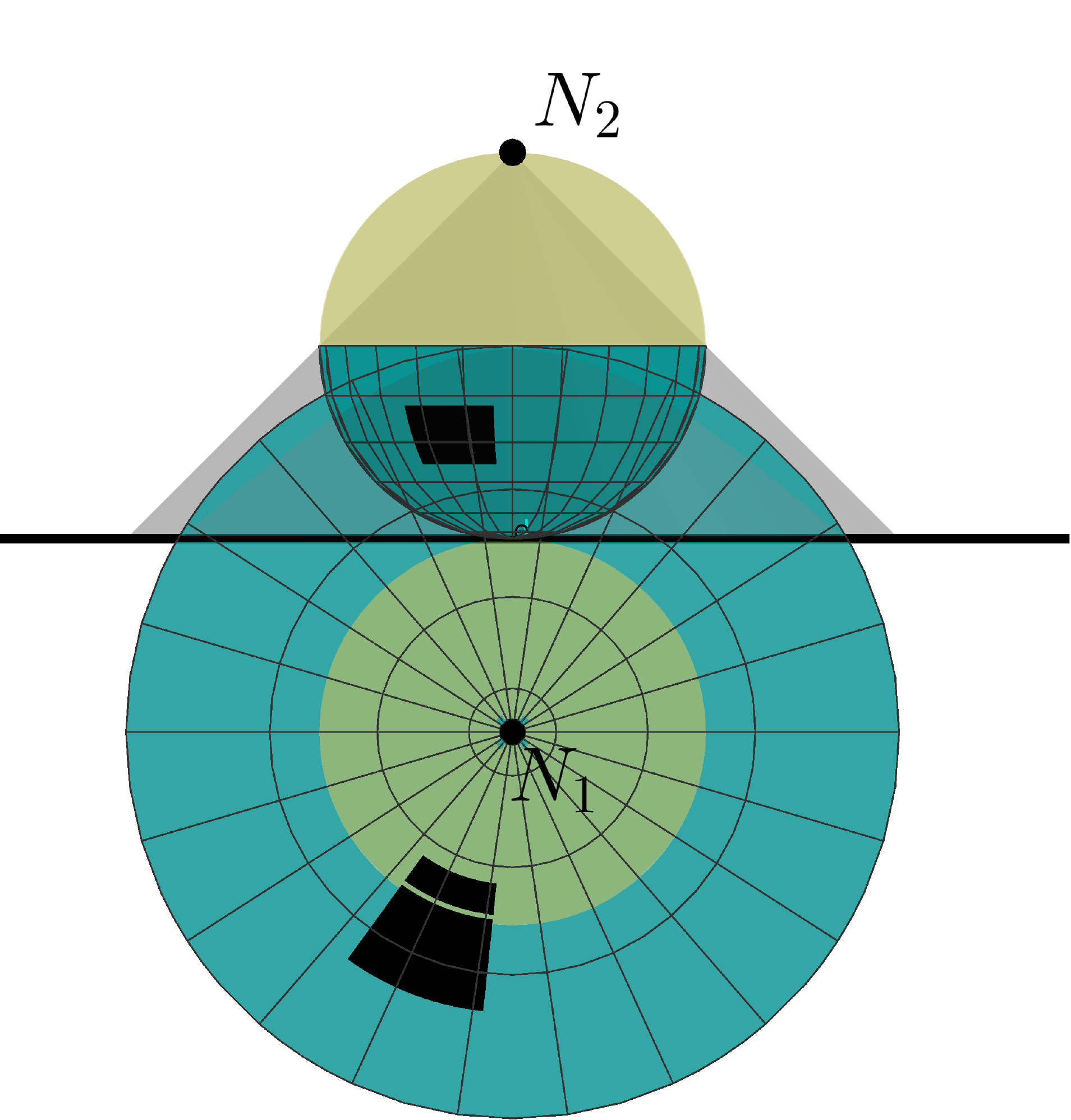}
\caption{A stereographic projection of the southern hemisphere from the North pole $N$ to the tangent plane at the South pole (i.e. the horizontal plane of projection), and its Monge's projection. Images of a quadrilateral with its sides along parallels and meridians are highlighted.\\
Interactive model available at \texttt{https://www.geogebra.org/m/mfctqfcs}}
\label{fig:sp3d}
\end{figure}

A projection between an $n$-sphere embedded in an $(n+1)$-space from its point $N$ and an $n$-dimensional hyperplane supplemented with a point at infinity $\{\infty\}$, not passing through $N$ is called stereographic projection. Due to its angle-preserving property, it is a convenient tool to visualize spheres. The 3-dimensional case of a stereographic projection of a 2-sphere from the North pole~$N$ to a tangent plane at the South pole and its Monge's projection is depicted in Figure~\ref{fig:sp3d}. For a brief overview see \cite{Snyder1987}, pp. 154--163, and for more general view with the construction of a spherical inversion used in Section~\ref{sec:sphinv} see \cite{Odehnal2020}, pp. 368--378.

A stereographic projection of a 3-sphere onto a 3-space was used to create 3-D printed models and study properties of four-dimensional polytopes in \cite{Schleimer2012} and later on in \cite{Hart2014} to visualize symmetries of the quaternionic group on a 3-D printed models of monkeys in a hypercube (limbs, head, and tail of the monkey represent faces of cubical cells). Well designed online applets with stereographic projections of a hypercube are at \cite{Sirius14000} and including other polytopes, sections, and tori in a 3-sphere at \cite{Cervone}. A stereographic projection of a 3-sphere is often used to visualize objects usually studied in topology. The topology of a 3-sphere is studied in \cite{Kocak1987}. Several animations and four-dimensional stereographic projection are described in \cite{Banchoff1996}, Chapters 6 and 8. The same author in \cite{Banchoff1988} used computer graphics to visualize stereographic images of Pinkall's tori corresponding in the Hopf fibration to simple closed curves on a 2-sphere (see also \cite{Pinkall1985} and \cite{3dxm}). An interactive visualization of a stereographic projection of a Clifford torus on a 3-sphere is at \cite{Balmens2012}. Videos and animations of the Hopf fibration and its stereographic projection are at \cite{Johnson2011, Dimensions, Chinyere2012}. Commented videos with interactive environment discussing quaternions and stereographic projection are at \cite{Eater}. A stereographic and double-stereographic projection of an arbitrary object in the 4-space was described in \cite{Ohori2017}. As an example of a recent application, a stereographic projection to a 3-sphere is used in \cite{Hosseinbor2015} to analyze multiple disconnected anatomical structures mathematically represented as a composition of compact finite three-dimensional surfaces.

\begin{figure}[!b]
\centering
\includegraphics[height=7cm, trim=50 50 50 50, clip]{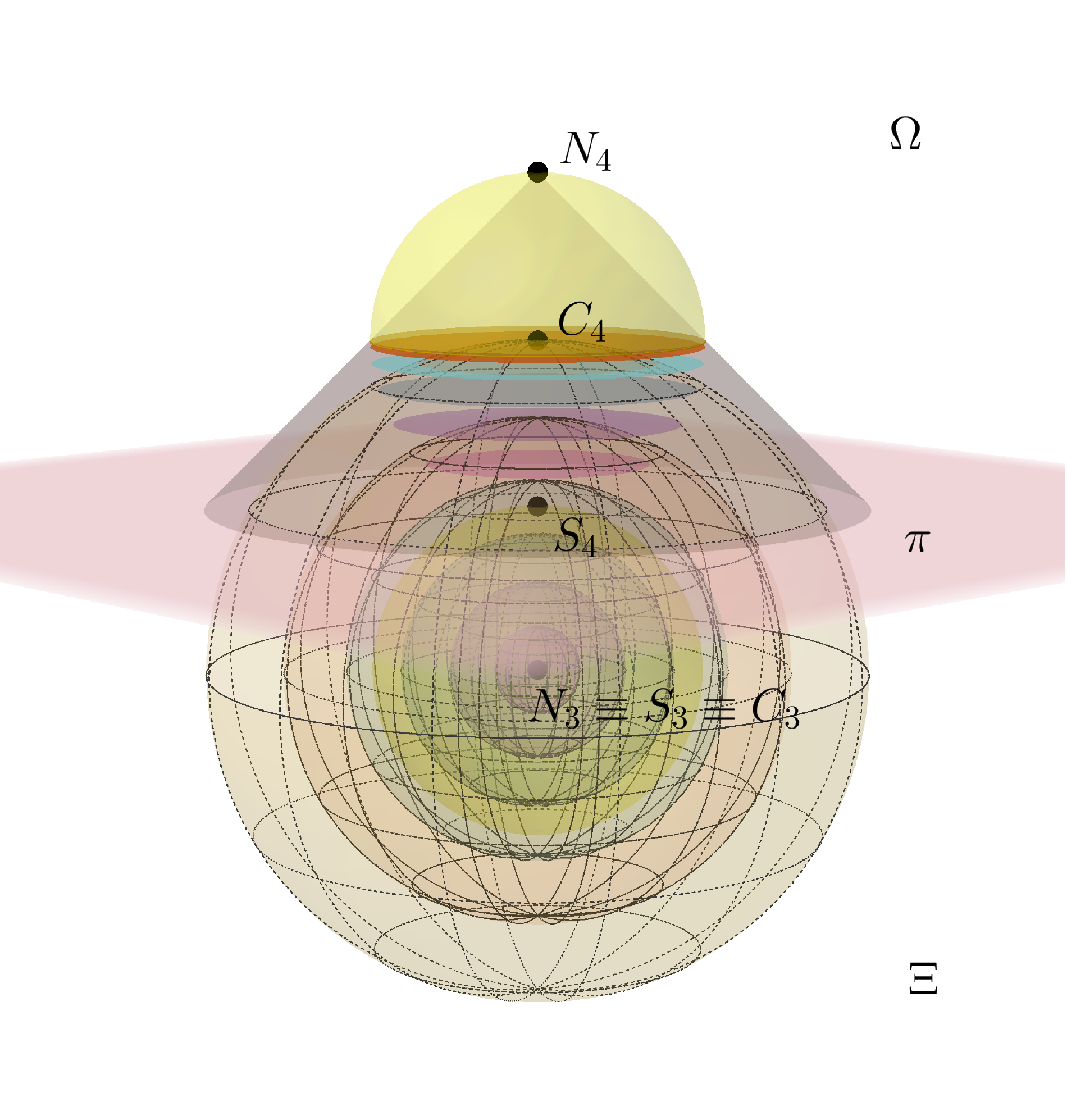}\hfill\includegraphics[height=7cm, trim=50 50 50 50, clip]{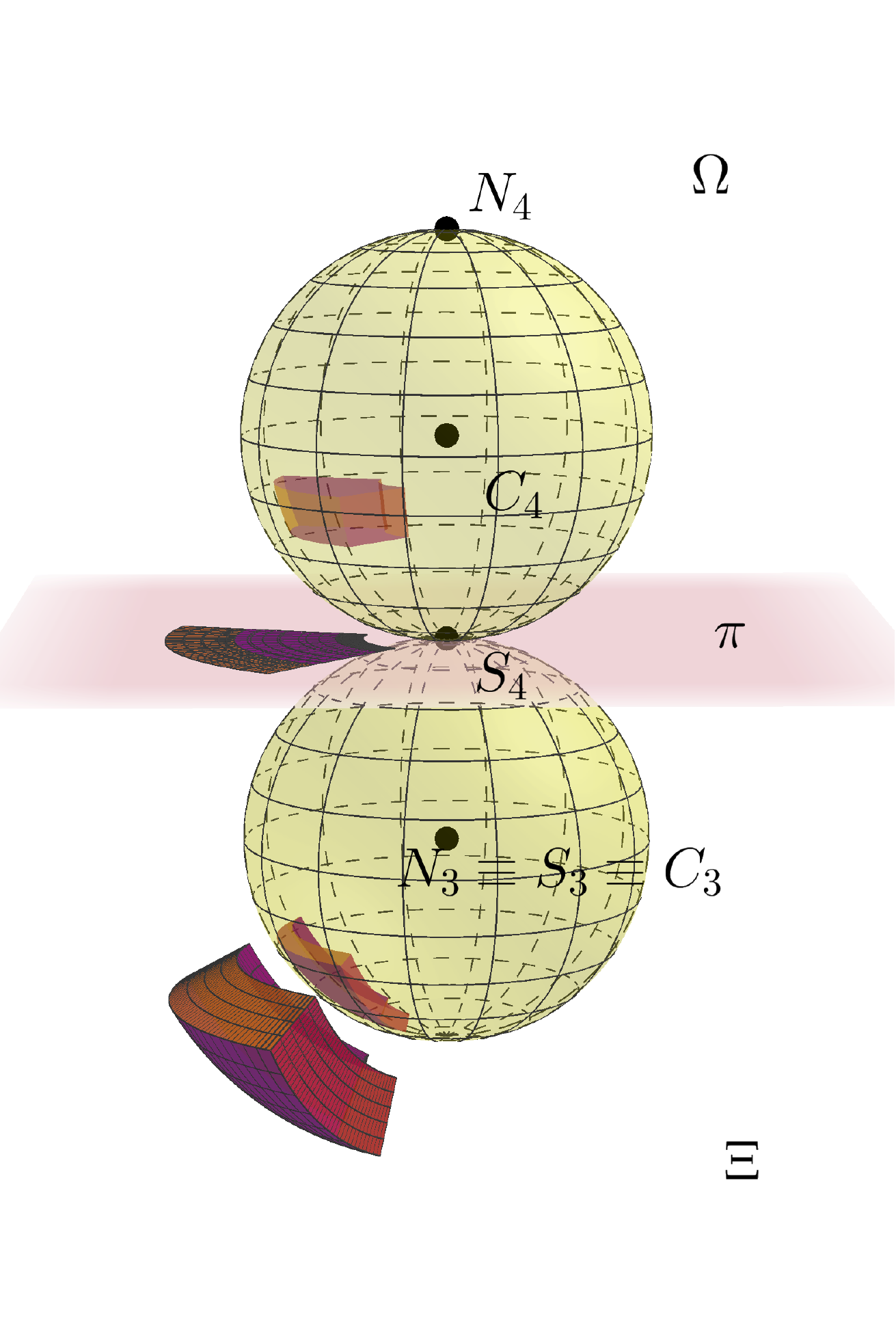}
\caption{The system of concentric spheres, which are stereographic images of parallel sections of the half of the 3-sphere split by the 3-space parallel to $\Xi$ from the point $N$ with the maximal $w$-coordinate to the tangent 3-space $\Xi$ at the antipodal point $S$, and its double orthogonal projection. On the right side are conjugated images and stereographic image of a hyperspherical hexahedron, as an analogy to a spherical quadrilateral.\\
Interactive model available at \texttt{https://www.geogebra.org/m/sz7ykk54}}
\label{fig:sp4d}
\end{figure}

The previous references were based on the analytic representation of points in the fourth dimension. With the use of the double orthogonal projections, we can construct images of four-dimensional objects synthetically and use the modeling 3-space to be also the projecting 3-space of a stereographic projection. A construction of a hyperspherical tetrahedron in a special position and its images in a stereographic and the double orthogonal projection is in \cite{Zamboj2020a}. A synthetic construction and animation of Hopf fibers in the double orthogonal projection and their stereographic projection are discussed in \cite{Zamboj2020b}. In this paper, we will extend these results in several aspects.

Let us describe the visualization of a 3-sphere and stereographic images in the double orthogonal projection in Figure~\ref{fig:sp4d} (left) in an analogy to the three-dimensional case in Figure~\ref{fig:sp3d} (right). In Monge's projection, the conjugated images of a \mbox{2-sphere} are disks, in the double orthogonal projection, the conjugated images of a 3-sphere are balls. In 3-D, sections of a \mbox{2-sphere} with planes parallel to the horizontal plane of projection (parallels) are circles, and their images are segments in the vertical plane and circles in the horizontal plane. The stereographic projection from the North pole $N$ to the tangent plane at the South pole $S$ projects these parallels to the system of concentric circles in the drawing plane. In 4-D, sections of a 3-sphere with 3-spaces parallel to the 3-space of projection $\Xi$ are 2-spheres, their $\Omega$-images are disks and $\Xi$-images are 2-spheres. The stereographic projection from the point $N$ with the maximal $w$-coordinate to the tangent 3-space $\Xi$ at the antipodal point $S$ projects the parallel 2-spheres to the system of concentric 2-spheres in the modeling 3-space. In both cases, the point $N$ is projected stereographically into the point at infinity $\{\infty\}$, and hence images of circles through $N$ become lines. In Figure~\ref{fig:sp4d} (right) (see also animated model) is a hyperspherical hexahedron along hyperspherical coordinates.

\section{A synthetic construction of a stereographic projection in the double orthogonal projection}

\begin{figure}[!htb]
\centering
\includegraphics[height=7cm, trim=100 50 50 30, clip]{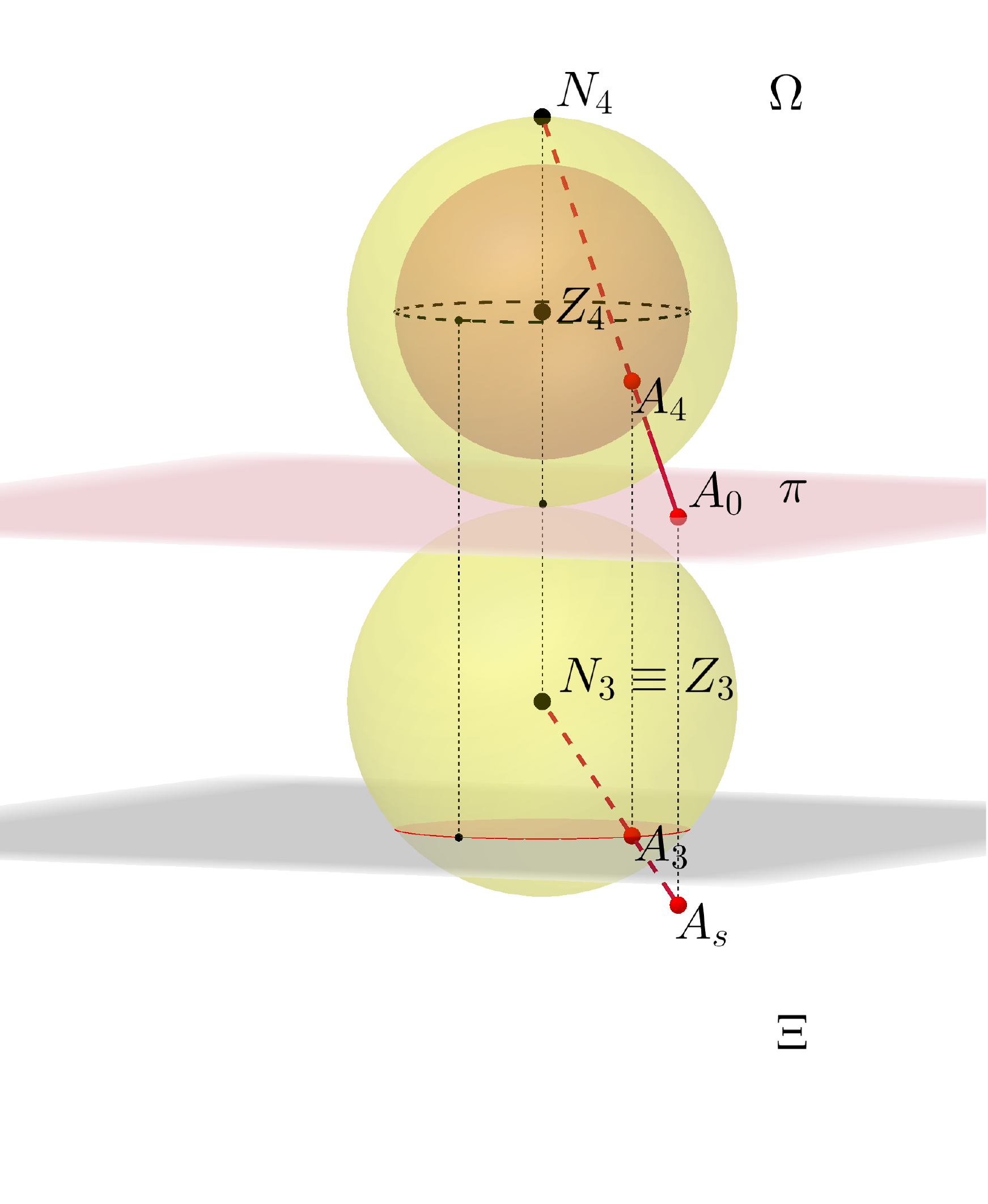}\hfill\includegraphics[height=7cm, trim=100 50 50 30, clip]{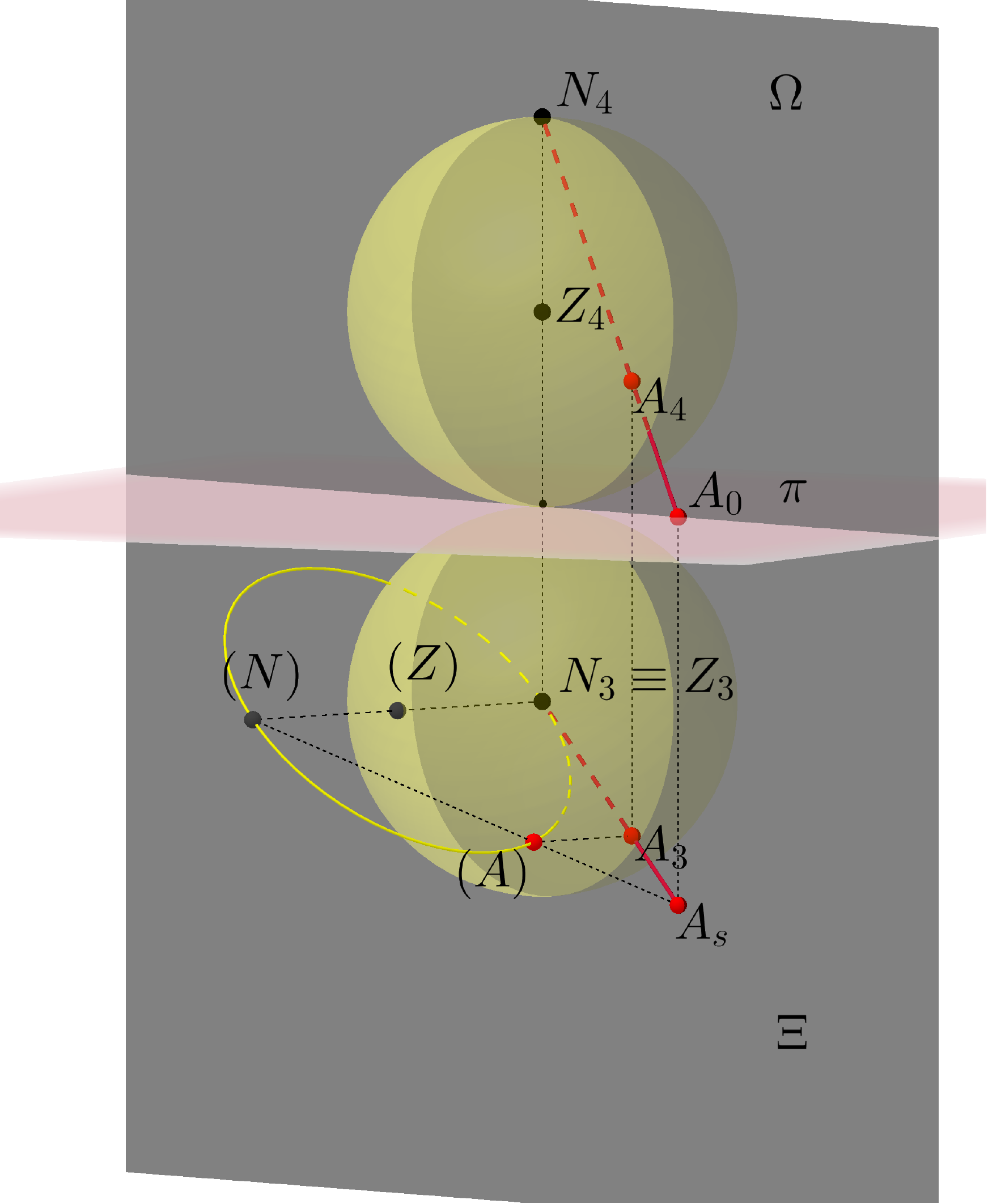}
\caption{(left) Construction of the stereographic image $A_s$ of the point $A$ on the 3-sphere given by its $\Omega$-image $A_4$. (right) Construction of the conjugated images $A_3$ and $A_4$ of a point $A$ from its stereographic image $A_s$.\\
Interactive model available at \texttt{https://www.geogebra.org/m/xdypddf9}}
\label{fig:sp4dpoint}
\end{figure}
The stereographic image of a point on a 3-sphere is the intersection of the projecting ray with the plane of projection (Figure~\ref{fig:sp4dpoint}, left). Let us have conjugated images of a 3-sphere with a center $Z$ in the double orthogonal projection and the $\Omega$-image $A_4$ of a point $A$ on the 3-sphere. The $\Xi$-image $A_3$ lies on the perpendicular to $\pi$ through $A_4$, i.e. ordinal line of the point $A$. Furthermore, the section of the 3-space through $A$ parallel to the 3-space $\Omega$ with a 3-sphere is a 2-sphere. Its $\Omega$-image is a 2-sphere with the center $Z_4$ through $A_4$, and its $\Xi$-image is a disk with the same radius in a plane parallel to $\pi$. Let $N$ be the abovementioned center of the sterographic projection to the 3-space $\Xi$. The stereographic image $A_s$ of the point $A$ on the 3-sphere lies on the line $N_3A_3$ and also on the perpendicular to $\pi$ through the intersection $A_0$ of $\pi$ and the line $N_4A_4$. Oppositely, the inverse construction of the conjugated images $A_3$ and $A_4$ from the stereographic image~$A_s$ is in Figure~\ref{fig:sp4dpoint} (right). In this case, we need to find the intersection~$A$ of the projecting ray $NA_s$ with the 3-sphere. For such construction, we can use a third orthogonal projection into the plane perpendicular to $\Xi$ through the projecting ray, and rotate it into the modeling 3-space (see \cite{Zamboj2019b} for more details). In the third view, the image of the 3-sphere is a circle with the center $(Z)$, and its intersection with the rotated line $(N)A_s$ is the rotated point $(A)$. Then, the point $A_3$ is constructed with the reverse rotation.

\begin{figure}[!htb]
\centering
\includegraphics[height=7cm, trim=100 50 50 30, clip]{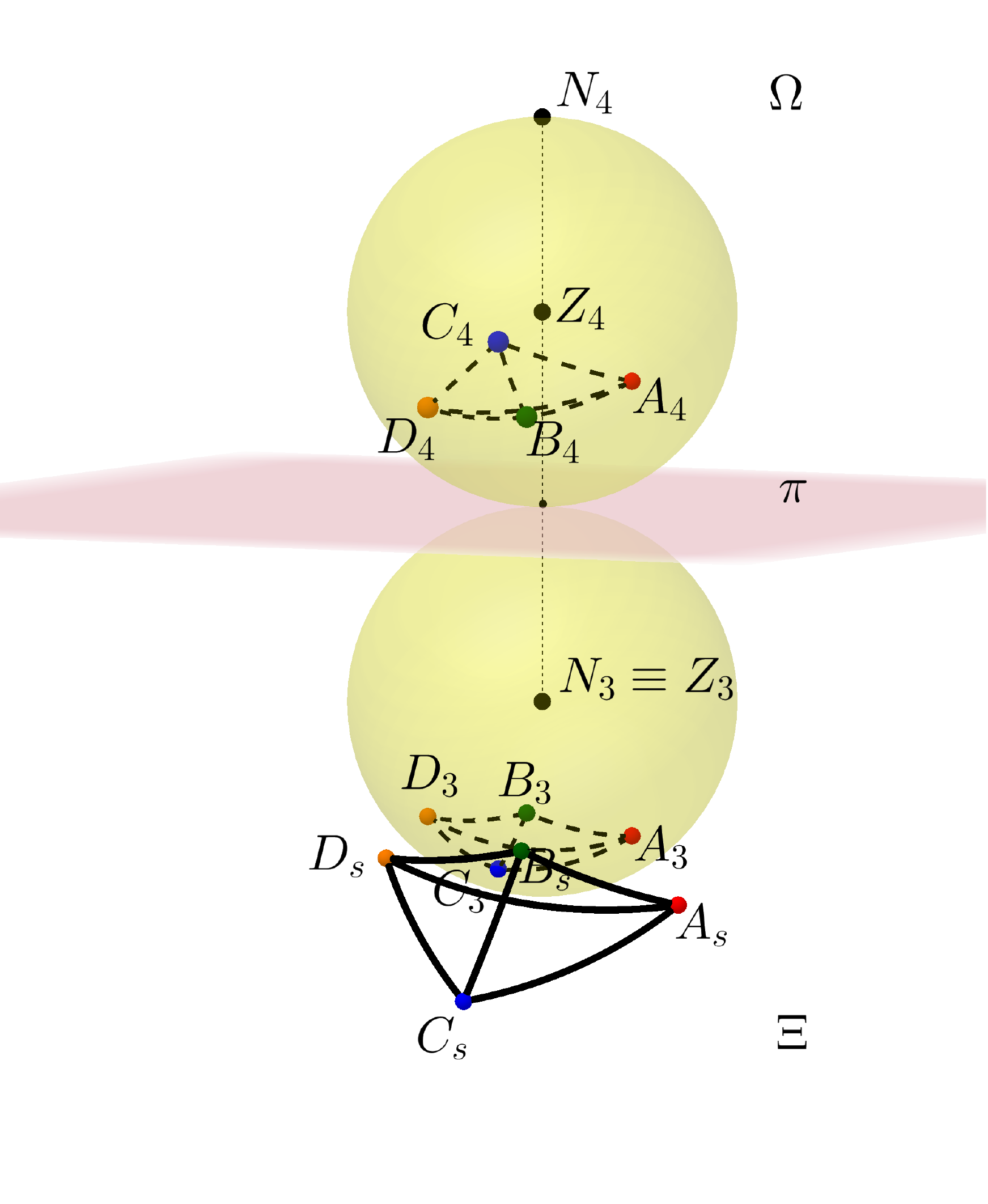}\hfill\includegraphics[height=7cm, trim=100 50 50 30, clip]{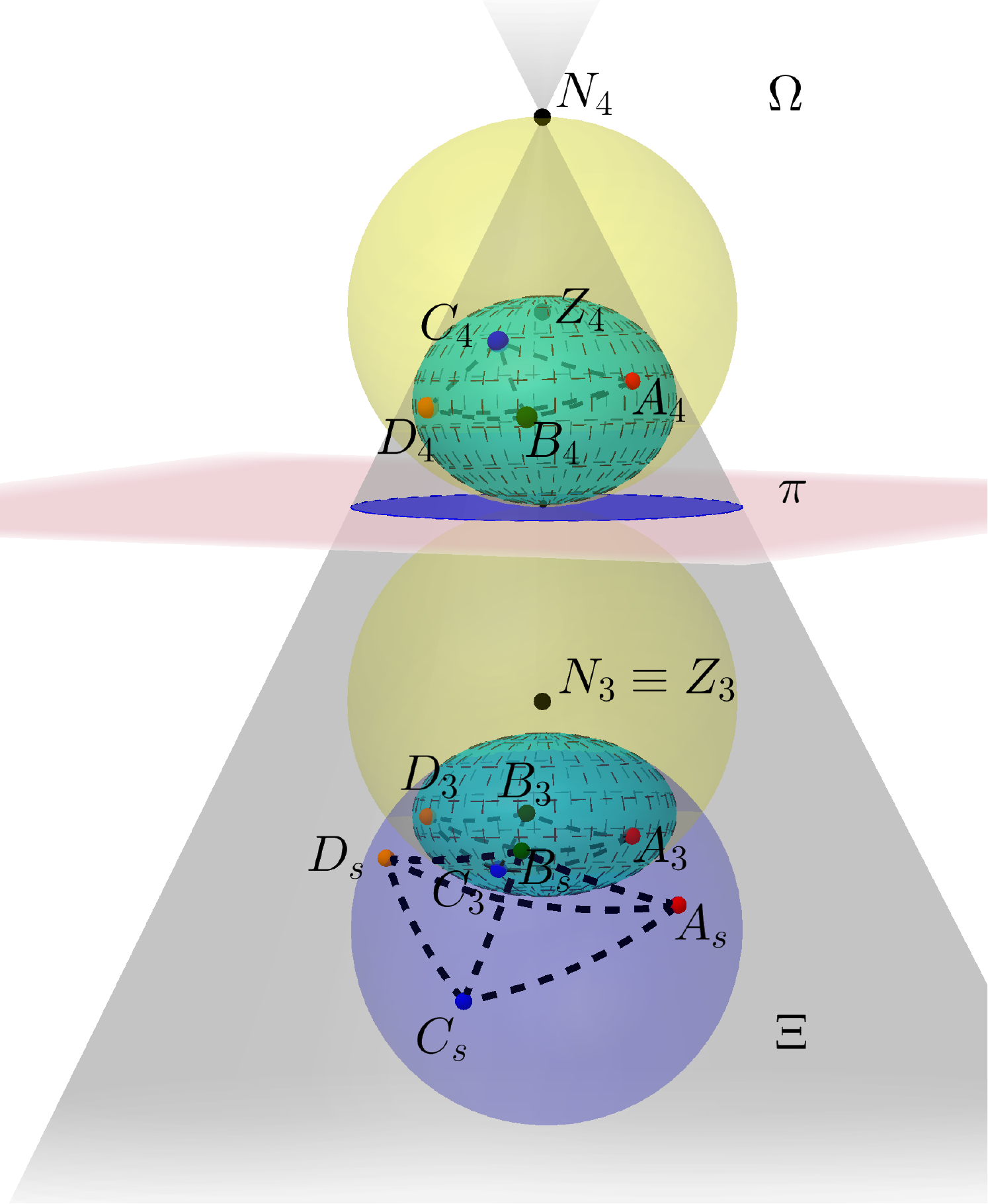}
\caption{(left) The double orthogonal projection and stereographic image of a hyperspherical tetrahedron $ABCD$.  (right) The same situation with the circumscribed 2-sphere around $ABCD$.\\
Interactive model available at \texttt{https://www.geogebra.org/m/xdypddf9}}
\label{fig:sp4dtetrahedron}
\end{figure}

In Figure~\ref{fig:sp4dtetrahedron} (left), the same method is used to construct the stereographic images $A_s, B_s, C_s, D_s$ of the vertices of a hyperspherical tetrahedron $ABCD$ (generalization of a spherical triangle). The edges of the tetrahedron $ABCD$ are circular arcs and faces are spherical triangles.\footnote{The faces are not depicted due to insufficient possibilities of the surface parametrization in GeoGebra, but the reader can turn on the visibility of the corresponding spheres in the stereographic projection in the online model.} These properties are preserved due to the conformity of stereographic projection. To construct the stereographic images of the circular edges, we can conveniently use stereographic images of the antipodal points. For example, the edge $A_sB_s$ lie on the circle $A_sB_sA'_s$, where $A'_s$ is the stereographic image of the point $A'$ antipodal to $A$, and so its conjugated images $A'_3$ and $A'_4$ are the mirror images of $A_3$ and $A_4$ about $Z_3$ and $Z_4$, respectively. The conjugated images of the edges are constructed point-by-point from their stereographic images and create elliptical arcs in a general position. Note that we could construct them, however, more laboriously, as the intersections of the 3-sphere and planes.

Further on, a 2-sphere circumscribed around a tetrahedron $ABCD$ is visualized in Figure~\ref{fig:sp4dtetrahedron} (right). While in the stereographic projection it is simply a 2-sphere in the true shape, the conjugated images in the double orthogonal projection are ellipsoids. They can be constructed as the intersections of a 3-sphere with the 3-space defined by the noncospatial points $A,B,C$, and $D$ (see \cite{Zamboj2019b} for details), or, as highlighted in the $\Omega$-image, as the intersections of the conical hypersurface with the vertex $N$ with a 3-sphere.

\begin{figure}[!htb]
\centering
\includegraphics[height=7cm]{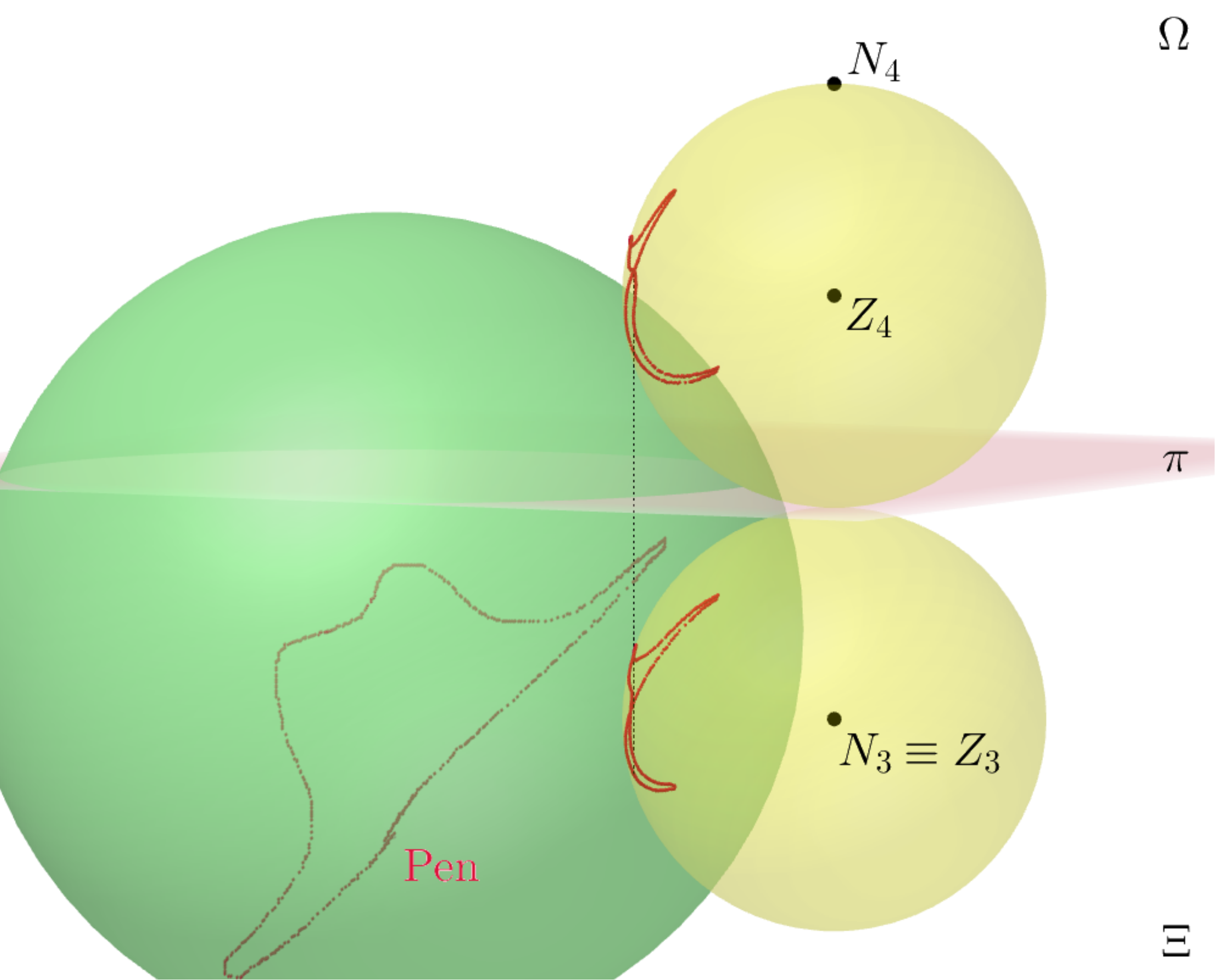}
\caption{My drawing was not a picture of a hat. It was a picture of a boa constrictor digesting an elephant on a 3-sphere in the double orthogonal projection (slightly modified from \cite{Exupery1943}).\\
Interactive model available at \texttt{https://www.geogebra.org/m/mjpeaud4}}
\label{fig:freehand}
\end{figure}

The interactive environment, in combination with the construction of a stereographic projection from a 3-space onto a 3-sphere (Figure~\ref{fig:sp4dpoint}, right), brings a possibility of freehand drawing on a 3-sphere. With the use of the trace tool in GeoGebra (or similar point-by-point construction in other software), we move the point in the stereographic projection, and its dependent conjugated images draw their traces. This way, we can draw a stereographic image (its $\Xi$-image) of a curve, or any other picture, in the modeling 3-space, and simultaneously its double orthogonal projection. A four-dimensional curve on a 3-sphere drawn by hand on a 2-sphere in a stereographic projection is in Figure~\ref{fig:freehand}. 

\subsection{Spherical inversion}
\label{sec:sphinv}

\begin{figure}[!htb]
\centering
\includegraphics[width=0.9\textwidth, trim= 40 0 0 0, clip]{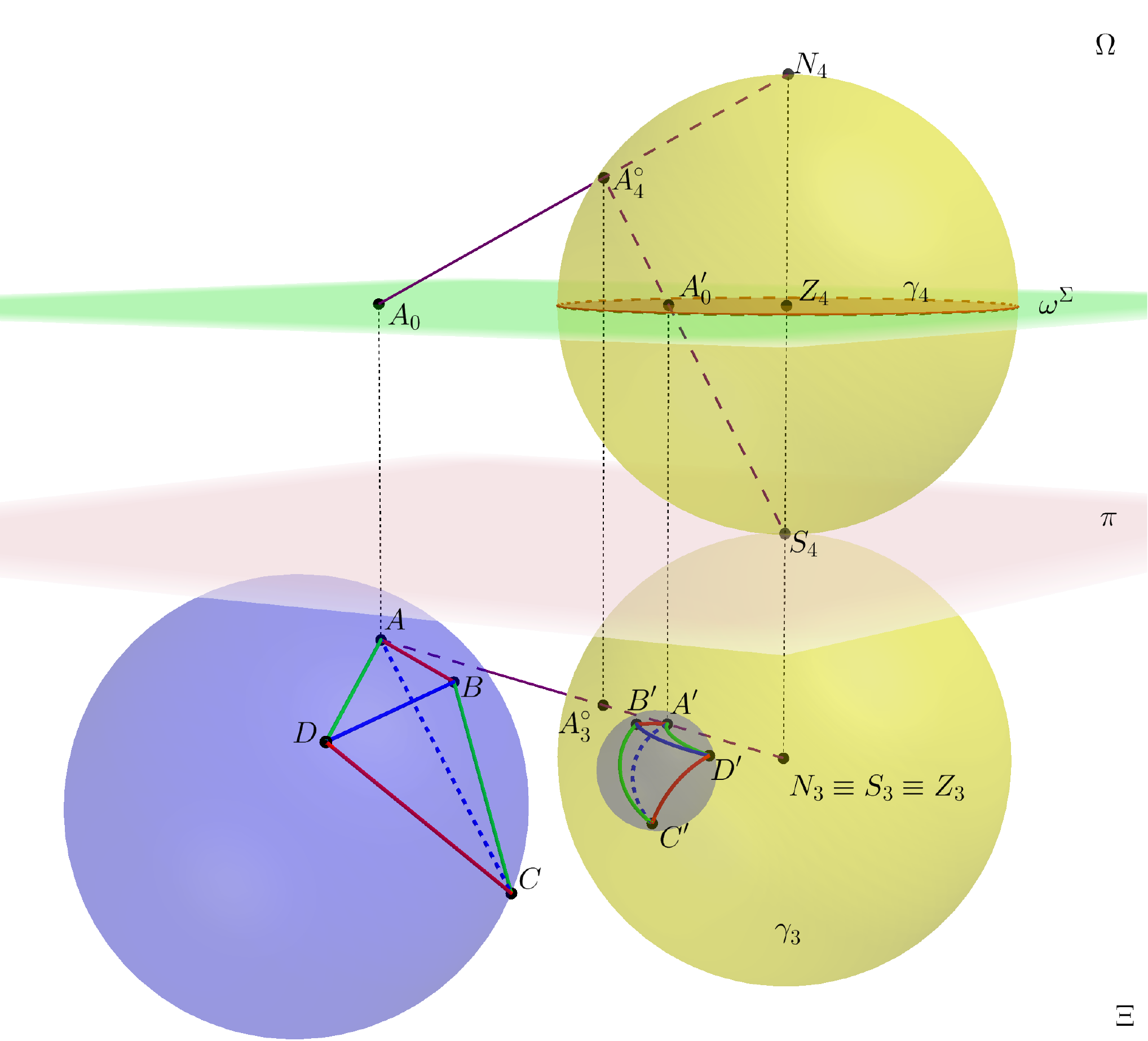}
\caption{Spherical inversion of a tetrahedron $ABCD$ and its circumscribed 2-sphere as a composition of stereographic projections from antipodal poles between a \mbox{3-sphere} and its equatorial 3-space.\\
Interactive model available at \texttt{https://www.geogebra.org/m/uqvek39f}}
\label{fig:sphinv}
\end{figure}

A spherical inversion may be obtained as a composition of two stereographic projections (Figure~\ref{fig:sphinv}). For this purpose, we choose the projecting 3-space $\Sigma$ of the stereographic projection through the center $Z$ of the 3-sphere and parallel to $\Xi$. A point~$A$ in the 3-space $\Sigma\cup\{\infty\}$ is stereographically projected in the projection~$p_N$ from the center $N$ to the point $A^\circ$ on the 3-sphere. The second stereographic projection $p_S$ from the center $S$ projects the point $A^\circ$ to the point $A'$ in the 3-space~\mbox{$\Sigma\cup\{\infty\}$}. Since $\Sigma$ is parallel to $\Xi$, all the $\Xi$-images of objects in the 3-space~$\Sigma$ are in the true shape. The section $\gamma$ of a 3-sphere and the 3-space~$\Sigma$ is a 2-sphere overlapping with the $\Xi$-image of the 3-sphere with the center $Z_3$ in the modeling 3-space. The composition of $p_N$ and $p_S$ is the spherical inversion, and $p_S(p_N(A))=A'$. Especially, all the points on $\gamma$ are fixed. Moreover, $p_S(p_N(\infty))=p_S(N)=Z$, and $p_S(p_N(Z))=p_S(S)=\infty$. 
\begin{figure}[!htb]
\centering
\includegraphics[height=7cm]{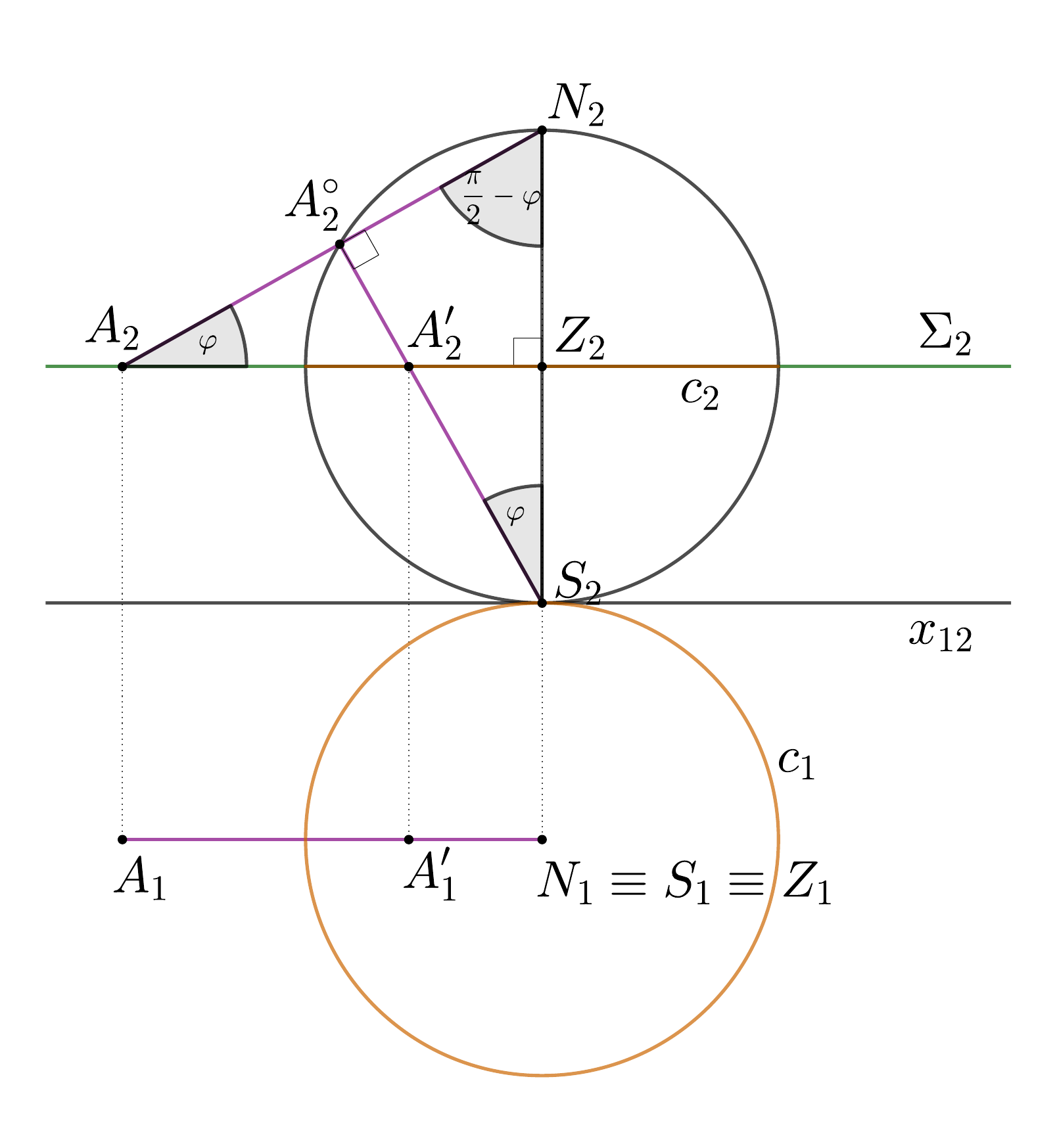}
\caption{Spherical inversion of a point $A$ in the orthogonal projection to the projecting plane of the line $\overline{NA}$.}
\label{fig:circinv}
\end{figure}
Note, that for each point $A\neq Z,\{\infty\}$ in $\Sigma$, we could choose a projecting plane of the line $AZ$ perpendicular to $\Xi$ which cuts the 2-sphere $\gamma$ in a circle $c$, and the final composition would be a circle inversion. The situation is depicted in Monge's projection in Figure~\ref{fig:circinv}, in which the circle inversion is in the orthogonal projection into the horizontal plane. Triangles $A_2Z_2N_2$ and $S_2Z_2A'_2$ are in their true shape in the front view, and they are apparently similar. Therefore, 
\begin{equation*}
\frac{|A_2Z_2|}{|Z_2N_2|}=\frac{|S_2Z_2|}{|Z_2A'_2|},
\end{equation*} and so 
\begin{equation*}
\label{eq:st}
|A_2Z_2||A'_2Z_2|=|Z_2N_2||Z_2S_2|=r^2,
\end{equation*} 
where $r$ is the radius of $c,\gamma$, and also the 3-sphere. From the given construction, it also holds that $|A_2Z_2|=|A_1Z_1|=|AZ|$ and $|A'_2Z_2|=|A'Z|$, and we have
\begin{equation*}
|AZ||A'Z|=r^2.
\end{equation*}
The last formula leads to the standard definition of the spherical inversion about the 2-sphere $\gamma$ with the center $Z$ and radius $r$.

Figure~\ref{fig:sphinv} also shows a tetrahedron $ABCD$ with a circumscribed sphere and their inversion. The conformity of stereographic projections is inherited in their composition to the spherical inversion. The lines not passing through the center of the inversion become circles. Consequently, the edges of the image tetrahedron $A'B'C'D'$ are circular arcs with preserved mutual angles. The image of the sphere circumscribed to $ABCD$ is a sphere circumscribed to $A'B'C'D'$.

\subsection{Hopf tori of Clelia curves}

\begin{figure}[!htb]
\centering
\includegraphics[width=0.31\textwidth]{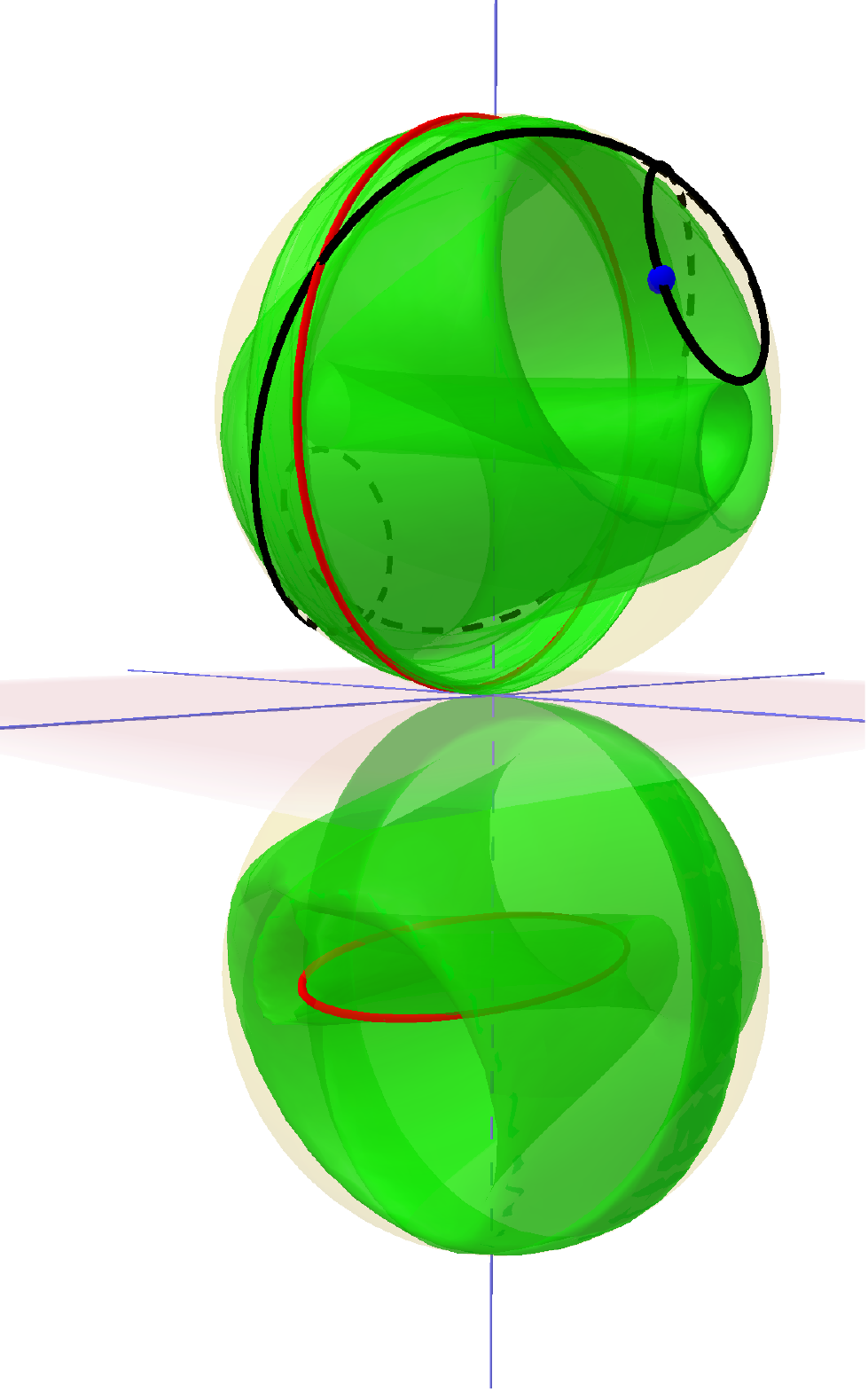} \hspace{1pt} \includegraphics[width=0.31\textwidth]{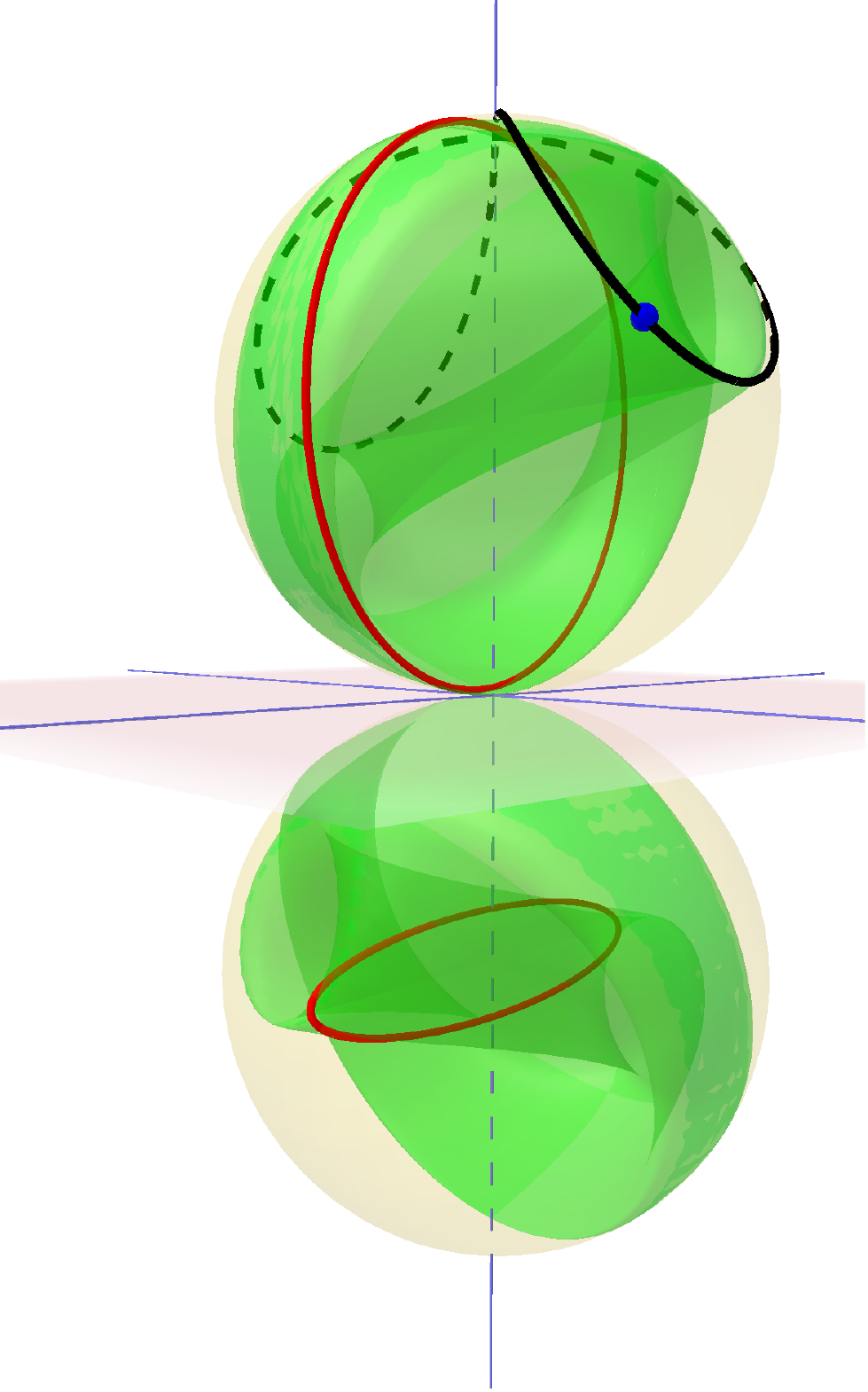} \hspace{1pt}\includegraphics[width=0.31\textwidth]{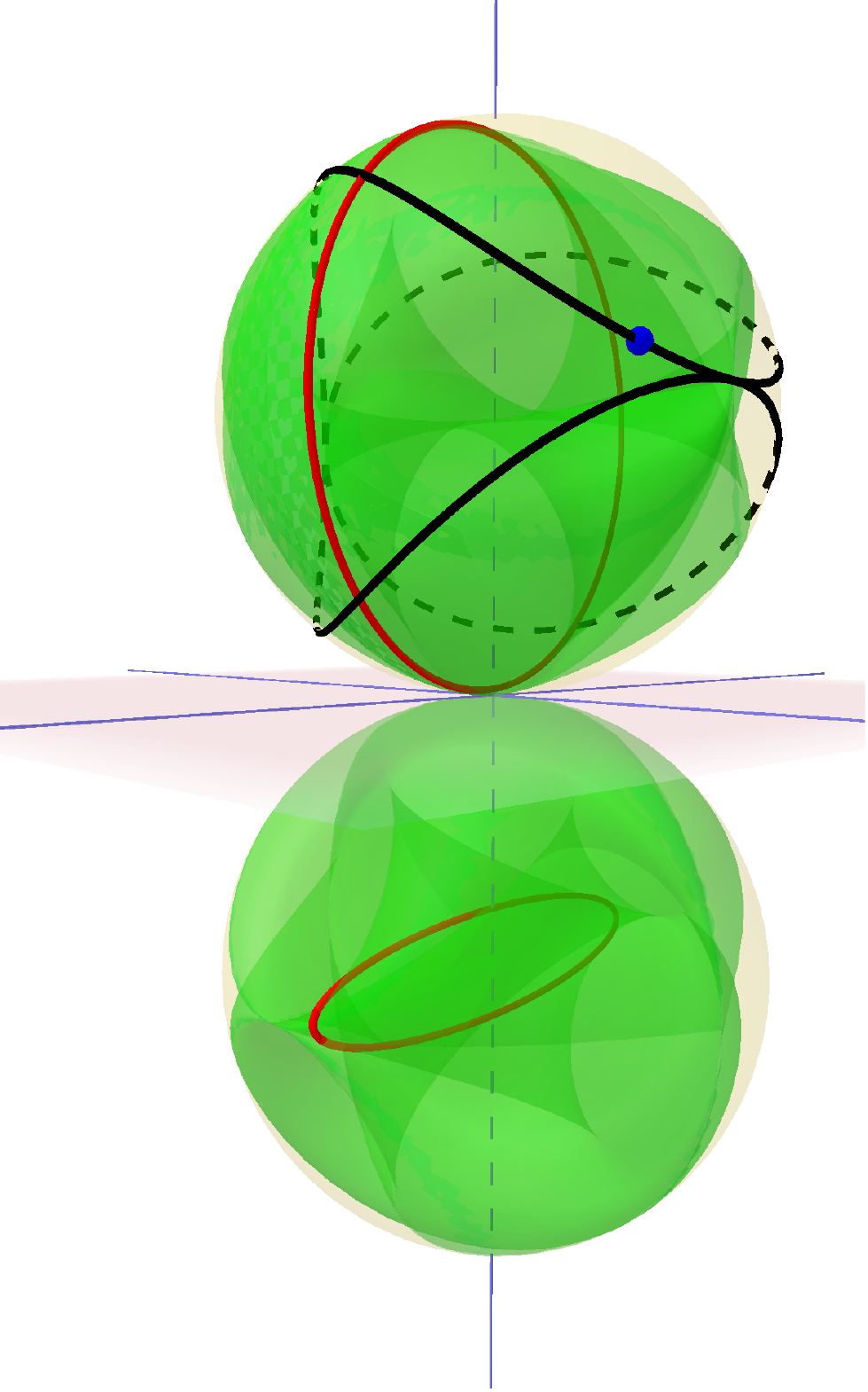}\\
\includegraphics[width=0.31\textwidth]{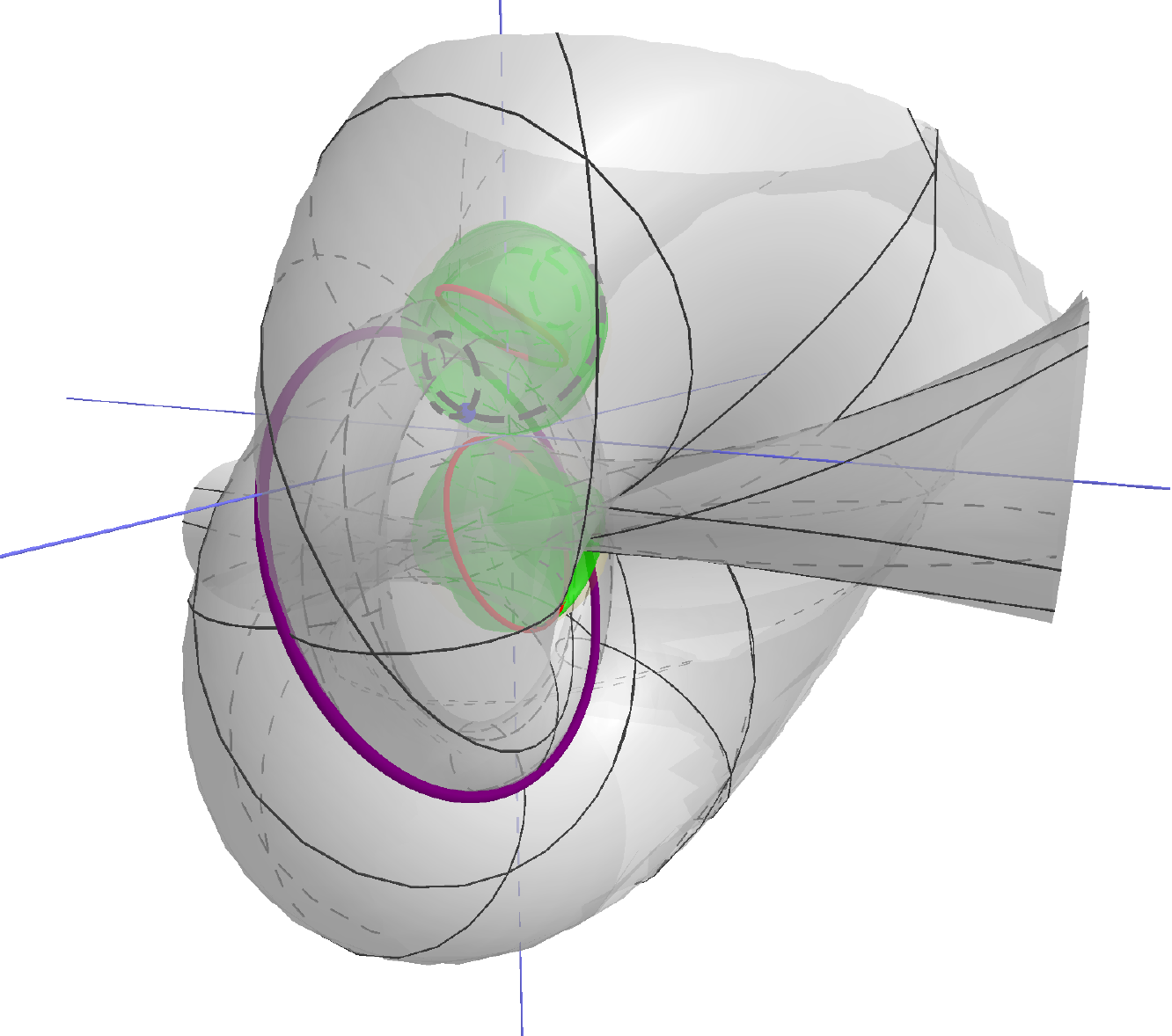} \hspace{1pt} \includegraphics[width=0.31\textwidth]{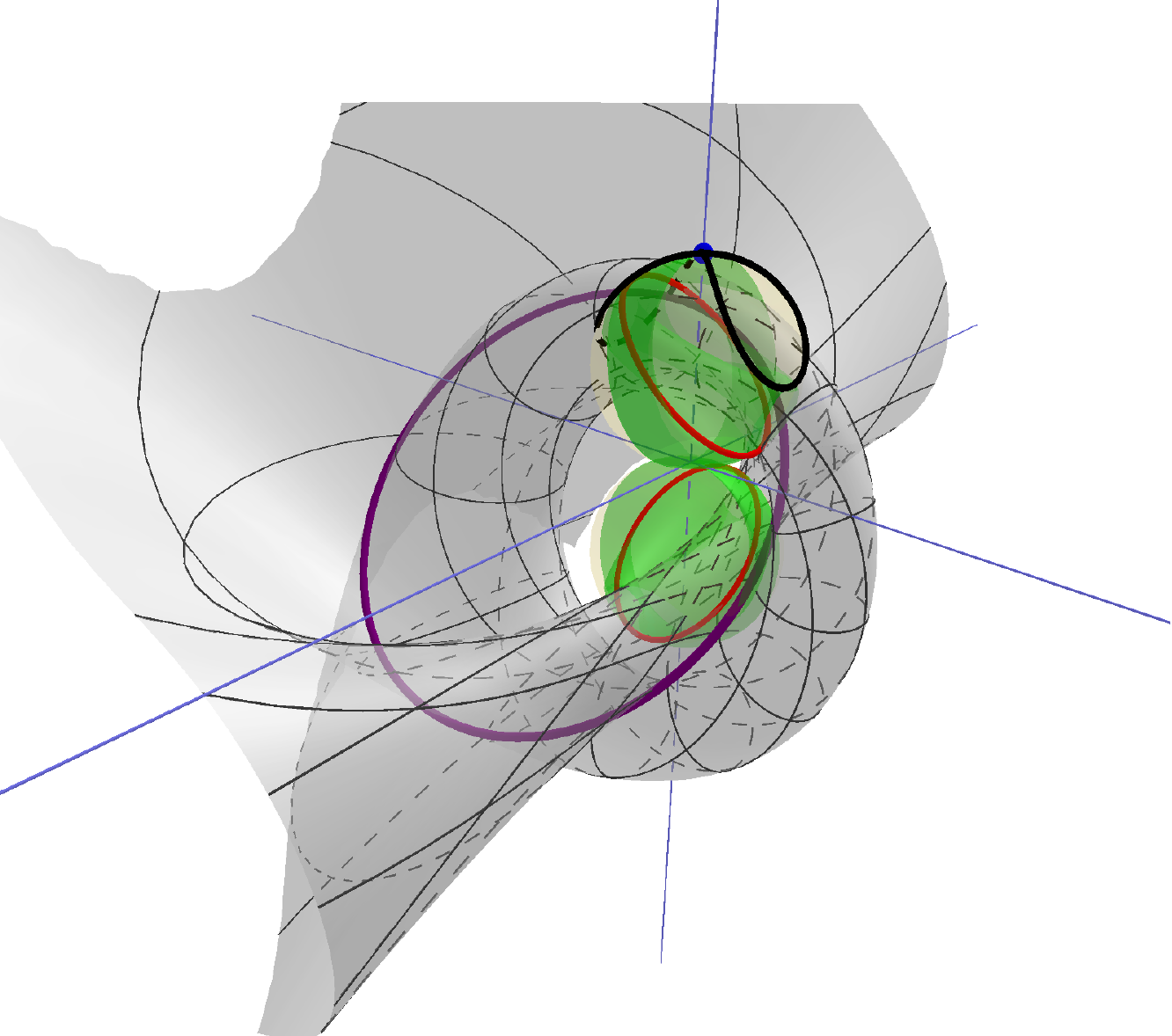} \hspace{1pt}\includegraphics[width=0.31\textwidth]{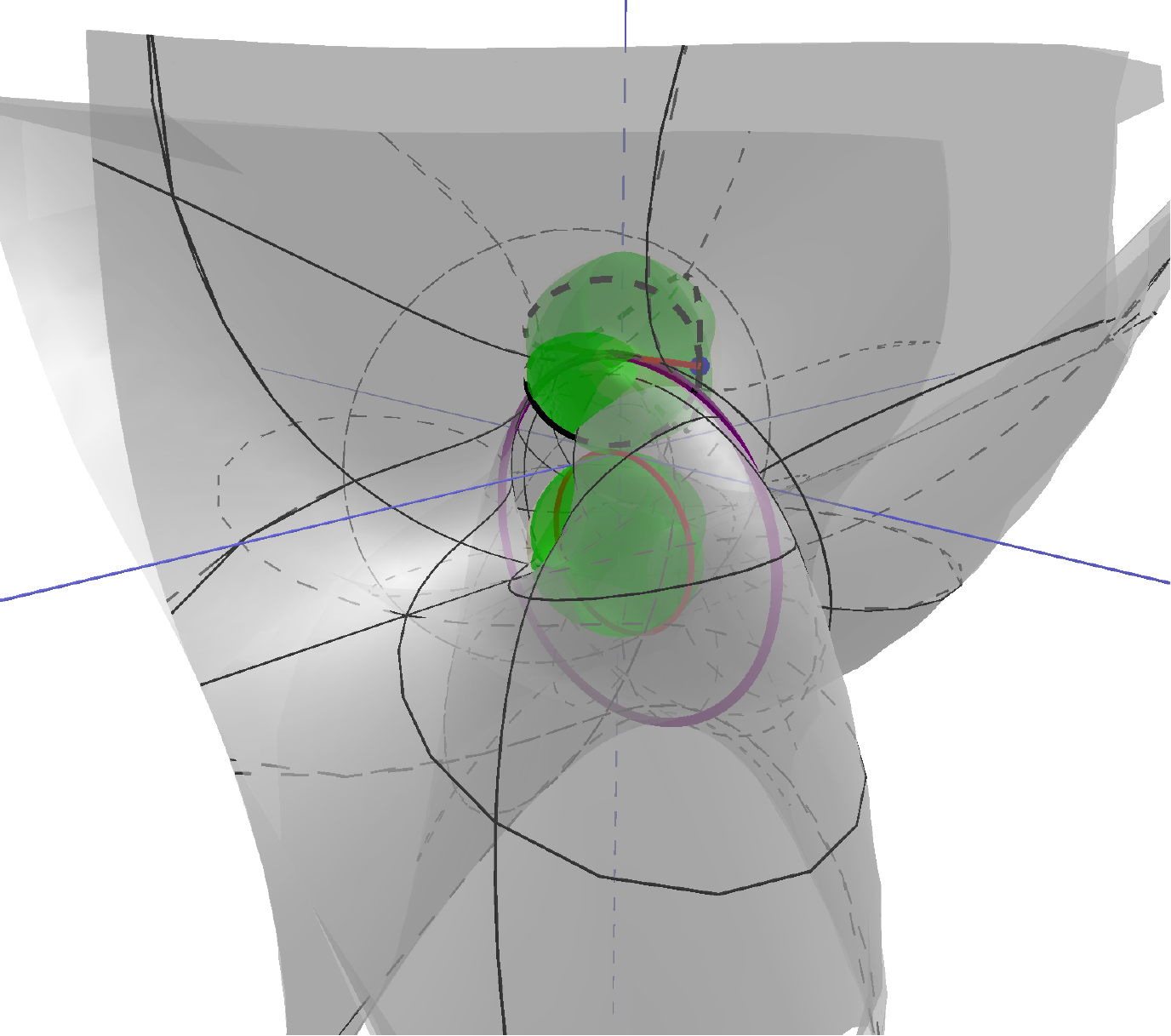}
\caption{(Up) Conjugated images of the Hopf tori generated by the Clelia curves for $s=0.5$ (Left), $s=1$ Viviani curve (Center), $s=2$ (Right). Highlighted is a point (blue) on the Clelia curve (black), conjugated images of its fiber (red) and conjugated images of the Hopf torus (green). (Down) Stereographic image of a fiber (purple) and Hopf torus (gray).\\
Interactive model available at \texttt{https://www.geogebra.org/m/sbf4xadp}}
\label{fig:clelia}
\end{figure}

The Hopf fibration is a mapping between a 3-sphere and a 2-sphere. In particular, it sends a point $$\displaystyle P(x_P,y_P,z_P,w_P)\in\mathbb{R}^4$$ on a 3-sphere to the point $$\displaystyle P'(2(x_Pz_P+y_Pw_P),2(-x_Pw_P+y_Pz_P),x_P^2-y_P^2-z_P^2-w_P^2)\in\mathbb{R}^3$$ on a 2-sphere. Oppositely, a point on a 2-sphere in spherical coordinates $$P'(\sin\psi_P \cos\varphi_P, \sin\psi_P\sin\varphi_P,\cos\psi_P)$$ for $\psi_P\in\langle0,\pi\rangle, \varphi_P\in\langle0,2\pi)$ corresponds to a set of points 
$$c_P(\cos\frac{\psi_P}{2}\cos(\varphi_P+\beta),\cos\frac{\psi_P}{2}\sin(\varphi_P+\beta),\sin\frac{\psi_P}{2}\cos\beta,\sin\frac{\psi_P}{2}\sin\beta)$$
for $\beta\in \langle0,2\pi)$, which forms a great circle (fiber) on a 3-sphere. If $P'$ lies on a closed curve on the 2-sphere, its corresponding Hopf fibers form a Hopf torus on the 3-sphere. For the visualization of Hopf tori, it is convenient to use the stereographic projection that preserves circles. This way we can construct and study tori on the 3-sphere in the 4-space given by a curve on the 2-sphere in the 3-space, and vice-versa.

For the sake of visualization in the double orthogonal projection, we perform several adjustments. First, we swap the reference axes $y$ and $z$, and hence the 3-spaces of projection will be $\Xi(x,y,z)$ (upper) and $\Omega(x,z,w)$ (lower)  with the common plane $\pi(x,z)$.\footnote{This choice reflects the possibility of a definition of the Hopf fibration in the complex number plane, and so the common plane $\pi(x,z)$ corresponds to the real parts of coordinates $(x+\i y, z+\i w)$ of points.} To avoid overlapping of the conjugated images, we translate the center of the abovementioned unit 3-sphere from the origin to the point $(0,1,0,1)$. The preimage 3-sphere and the image 2-sphere are not strictly related in the definition of the Hopf fibration, and topologists often describe them as separated objects. However, in our interpretation, the 2-sphere is chosen to be the equatorial section of the 3-sphere with the 3-space parallel to $\Omega(x,z,w)$, so its $\Xi$-image is in the true shape. Consequently, points of the 3-sphere are stereographically projected from the center $(0,2,0,1)$ to the 3-space $\Omega(x,z,w)$. A synthetic construction of the Hopf fiber of any point on the 2-sphere and a Hopf torus of a circle on a 2-sphere and their stereographic images in this setting are described in \cite{Zamboj2020b}. 

We extend these results and highlight interactive possibilities for the so-called Clelia curves given parametrically on the (translated) 2-sphere in the form $$k(\psi)=(\sin(s\psi)\cos\psi, \sin(s\psi)\sin\psi+1,\cos(s\psi)),$$ where $s\in\mathbb{R}$ defines the specific curve, for $\psi\in I\subset\mathbb{R}$ (see \cite{Goemans2016} for details and further generalizations). Conjugated images and stereographic projections of the Hopf tori corresponding to the Clelia curves for $s=\frac{1}{2}, 1,$ and $2$ are visualized in Figure~\ref{fig:clelia}. In the interactive model, the user can choose the parameter $s$ of the curve and change the corresponding Hopf torus, and also move with a point~$P$ on $k(\psi)$ and its corresponding fiber on the 3-sphere by manipulating the parameter~$\psi$. An example of a simple straightforward observation is, that a point, in which the curve on the 2-sphere intersects itself, becomes a circle (or line) in which the torus on the 3-sphere intersects itself. 

\section{Conclusion}
The combination of the double orthogonal projection and stereographic projection is an accessible tool for the investigation of a 3-sphere in the 3-space. The simplicity of the synthetic construction of the stereographic image of a point and the use of interactive 3-D modeling software imply that we can actually draw sketches and shapes on the 3-sphere embedded in the 4-space. We have also shown a construction of the well-known relationship between the stereographic projection and the spherical inversion. Our results were applied to the construction of Hopf tori on a 3-sphere generated by Clelia curves on a 2-sphere. The visualizations were presented on interactive 3-D models, which are easily extendible for further theoretical and practical applications in various software. 

%
%
\bibliographystyle{spmpsci}
\bibliography{127-paper-bibliography-Zamboj}








\end{document}